\documentclass{article}
\begin{document}
\newtheorem{proposition}{Proposition}[section]
\newtheorem{definition}{Definition}[section]
\newtheorem{lemma}{Lemma}[section]
 
\title{\bf Number-Like Objects and the Extended Lie Correspondence}
\author{Keqin Liu\\Department of Mathematics\\The University of British Columbia\\Vancouver, BC\\
Canada, V6T 1Z2}
\date{December, 2005---April, 2006}
\maketitle

\begin{abstract} We classify finite dimensional division real associative $\mathcal{Z}_2$-algebras, introduce composition $\mathcal{Z}_2$-algebras, and extend the Campbell-Baker-Hausdorff series and Lie correspondence in the context of linear Hu-Liu Leibniz algebras.
\end{abstract}

The impetus to my research is the curiosity to generalize the beautiful theory of Lie algebras. A natural criterion of judging whether a new algebraic object is a good generalization of a Lie algebra is whether there is a fair counterpart of the Lie correspondence between connected linear groups and linear Lie algebras. One way of introducing new algebraic objects which satisfy the criterion is to add some new algebraic structures to Lie algebras. In particular, if a Leibniz algebra structure is added to a nonsemisimple Lie algebra, then the resulting algebraic object can be used to establish the extended Lie correspondence. This is one of the results of this paper. Another result of this paper is the extension of Frobeius' theorem, which is an unexpected result appearing in my search for the generalization of Lie algebras. It is well known that Frobeius' theorem asserts that only finite dimensional division real associative algebras are the field $\mathcal{R}$ of the real numbers, the field $\mathcal{C}$ of the complex numbers, and Hamilton's quaternion algebra $\mathcal{H}$. In the context of finite dimensional division real associative $\mathcal{Z}_2$-algebras, the list in Frobeius' theorem is extended to a new list consisting of eight objects, where the first three objects are the ordinary finite dimensional division real associative algebras (without divisors of zero), and the last five objects are finite dimensional real associative $\mathcal{Z}_2$-algebras with divisors of zero. Hence, it is reasonable to regard finite dimensional division real associative $\mathcal{Z}_2$-algebras as number-like objects.

\medskip
Throughout this paper the associative algebras considered are assumed to have an identity.

\bigskip
\section{Basic Definisions}

In this section we recall some basic definitions introduced in \cite{hl1}. 

\begin{definition}\label{def1.1} A Lie algebra $( L,\, +,\, [ \,, \,])$ is called a {\bf Hu-Liu Leibniz algebra} if there exists a binary operation 
$\langle \,, \,\rangle : L\times L \to L$ such that the following two properties hold:
\begin{description}
\item[(i)] $\langle \,, \,\rangle$ satisfies the {\bf Leibniz identity}:
\begin{equation}\label{eq1}
\langle\langle x, y\rangle , z\rangle=\langle x, \langle y, z\rangle\rangle+ 
\langle\langle x,  z\rangle , y\rangle \qquad\mbox{for $x,y,z\in\mathcal{L}$.}
\end{equation}
\item[(ii)] $\langle \,, \,\rangle $ and $ [\, ,\, ]$ satisfy the following 
{\bf Hu-Liu identities}:
\begin{equation}\label{eq2}
\langle x, [y, z]\rangle =\langle x, \langle y, z\rangle\rangle,
\end{equation}
\begin{equation}\label{eq3}
[\langle x, x\rangle , y] =\langle \langle x, x\rangle , y\rangle,
\end{equation}
\begin{equation}\label{eq4}
\langle [x, y], z\rangle +[\langle y, z\rangle , x]+[y, \langle x, z\rangle ]=0,
\end{equation}
\begin{equation}\label{eq5}
[\langle x, y\rangle , z] +[z, [x, y]]+[z, \langle y, x\rangle  ]+
\langle z, \langle x, y\rangle\rangle=0
\end{equation}
for $x$, $y$, $z\in L$.
\end{description}
\end{definition}

A Hu-Liu Leibniz algebra $L$ is denoted by $( L,\, +,\, \langle \,, \,\rangle ,\, [\, ,\, ])$, where $\langle \,, \,\rangle $ and $ [\, ,\, ]$ are called the {\bf angle bracket} and the {\bf square bracket}, respectively. A subspace $I$ of a Hu-Liu Leibniz algebra 
$( L,\, +,\, \langle \,, \,\rangle ,\, [\, ,\, ])$ is called a {\bf Hu-Liu Leibniz subalgebra} of $L$ if $\langle I, I\rangle\subseteq I$ and $[I, I]\subseteq I$.

\begin{definition}\label{def1.2} Let $A=A_0\oplus A_1$ (as vector spaces) be an algebra.
\begin{description}
\item[(i)] $A$ is called a {\bf $\mathcal{Z}_2$-algebra} if 
\begin{equation}\label{eq6}
A_0A_0\subseteq A_0,\quad A_0A_1+A_1A_0\subseteq A_1 \quad\mbox{and}\quad A_1A_1=0.
\end{equation}
\item[(ii)] $A$ is called an {\bf associative $\mathcal{Z}_2$-algebra} if $A$ is both a $\mathcal{Z}_2$-algebra and an associative algebra.
\item[(iii)] $A$ is called a {\bf Banach $\mathcal{Z}_2$-algebra} if $A$ is both a associative $\mathcal{Z}_2$-algebra and a Banach algebra.
\end{description}
\end{definition}

Let $A=A_0\oplus A_1$ be a $\mathcal{Z}_2$-algebra, where the subspaces $A_0$ and $A_1$ satisfy (\ref{eq6}). $A_0$ and $A_1$ are called the {\bf even part} and {\bf odd part} of $A$, respectively. An element $a$ of $A$ can be written uniquely as $a=a_0+a_1$, where $a_0\in A_0$ and $a_1\in A_1$ are called the {\bf even component} and {\bf odd component} of $a$, respectively.

By Proposition 1.5 in \cite{Liu3}, if $A=A_0\oplus A_1$ is an associative $\mathcal{Z}_2$-algebra, then $( A,\, +,\, \langle \,, \,\rangle , \,[\, , \,])$ is a Hu-Liu Leibniz algebra, where the angle bracket $\langle \, , \, \rangle$ and the square bracket $[\, , \,]$ are defined by
\begin{equation}\label{eq7}
\langle x, y\rangle: =xy_0-y_0x \quad\mbox{and}\quad [x, y]:=xy-yx
\end{equation}
for $x\in A$, $y=y_0+y_1\in A$ and $y_i\in A_i$ with $i=0$ and $1$. A Hu-Liu Leibniz subalgebra of $( A,\, +,\, \langle \,, \,\rangle , \,[\, , \,])$ is called a {\bf linear Hu-Liu Leibniz algebra} in $A$. Let
$$
A^{-1}:=\{\, x\in A \, |\, \mbox{there exists $y\in A$ such that $yx=1=xy$ }\,\}
$$
and
$$
A_0^{-1}:=\{\, x_0\in A_0 \, |\, \mbox{there exists $y_0\in A_0$ such that $y_0x_0=1=x_0y_0$} \,\}.
$$
It is clear that
$$
A^{-1}=\{\, x_0+x_1 \, |\, \mbox{$x_0\in A_0^{-1}$ and $x_1\in A_1$ }\,\}.
$$
A subgroup $G$ of $A^{-1}$ is called a {\bf linear $\xi$-group} in $A$ if
\begin{equation}\label{eq8}
x_0\, G \, x_0^{-1}\subseteq G \quad
\mbox{for $x_0+x_1\in G$ with $x_0\in A_0$ and $x_1\in A_1$.}
\end{equation}

\bigskip
\section{Composition $\mathcal{Z}_2$-Algebras}

Except the dual real numbers $\mathcal{R}^{(2)}$ in \cite{B}, there are four more division real associative $\mathcal{Z}_2$-algebras: $\mathcal{C}^{(2)}$, $\mathcal{C}^{(-2)}$, $\mathcal{H}^{(2)}$ and $\mathcal{H}^{(-2)}$. See \cite{Liu3} for their definitions. It turns out that these five division associative $\mathcal{Z}_2$-algebras are all finite dimensional division real associative $\mathcal{Z}_2$-algebras with nonzero odd parts, which is a corollary of the following extension of Frobeius' theorem.

\begin{proposition}\label{pr2.1} The only finite dimensional division associative $\mathcal{Z}_2$-algebras over the field of real numbers are:
$$\mbox{
(1) $\mathcal{R}$, (2) $\mathcal{C}$, (3) $\mathcal{H}$, (4) $\mathcal{R}^{(2)}$,
(5) $\mathcal{C}^{(2)}$, (6) $\mathcal{C}^{(-2)}$, (7) $\mathcal{H}^{(2)}$, and
(8) $\mathcal{H}^{(-2)}$.}$$
\end{proposition}

\hfill\raisebox{1mm}{\framebox[2mm]{}}

\bigskip
We now introduce the notion of a composition $\mathcal{Z}_2$-algebras.

\begin{definition}\label{def2.1} A real $\mathcal{Z}_2$-algebra $A=A_0\oplus A_1$ is called a {\bf composition $\mathcal{Z}_2$-algebra} if $A$ is a Banach space with a norm $||\quad ||$ such that 
\begin{eqnarray*}
 ||xy||&\le& ||x||\,||y||,\\
 ||x_0y_0||&=& ||x_0||\,||y_0||,\\
 ||x_0y_1||&=& ||y_1x_0||=||x_0||\,||y_1||,
\end{eqnarray*}
where $x$, $y\in A$, $x_0$, $y_0\in A_0$ and $y_1\in A_1$.
\end{definition}

Each of the eight finite dimensional division real associative $\mathcal{Z}_2$-algebras is a 
composition $\mathcal{Z}_2$-algebra. The next example gives two non-associative composition $\mathcal{Z}_2$-algebras.

\bigskip
\noindent
{\bf Example } Let $\mathcal{O}^{(2\lambda)}=\mathcal{O}^{(2\lambda)}_0\oplus \mathcal{O}^{(2\lambda)}_1$ be a 16-dimensional $\mathcal{Z}_2$-graded vector space with $\lambda=\pm 1$, where the even part $\mathcal{O}^{(2\lambda)}_0$ and the odd part $\mathcal{O}^{(2\lambda)}_1$ are given by 
$\mathcal{O}^{(2\lambda)}_i:=\displaystyle\bigoplus_{j=1}^{8}\mathcal{R}e_{ij}$ for $i=0$ and $1$. $\mathcal{O}^{(2\lambda)}$ becomes non-associative composition $\mathcal{Z}_2$-algebra under the mutiplication tables:
$$
\begin{tabular}{|r||r|r|r|r|r|r|r|r|}\hline
$\cdot$ &$e_{01}$&$e_{02}$&$e_{03}$&$e_{04}$&$e_{05}$&$e_{06}$&$e_{07}$&$e_{08}$\\
\hline\hline
$e_{01}$&$e_{01}$&$e_{02}$&$e_{03}$&$e_{04}$&$e_{05}$&$e_{06}$&$e_{07}$&$e_{08}$\\
\hline
$e_{02}$&$e_{02}$&$-e_{01}$&$e_{04}$&$-e_{03}$&$e_{06}$&$-e_{05}$&$-e_{08}$&$e_{07}$\\
\hline
$e_{03}$&$e_{03}$&$-e_{04}$&$-e_{01}$&$e_{02}$&$e_{07}$&$e_{08}$&$-e_{05}$&$-e_{06}$\\
\hline
$e_{04}$&$e_{04}$&$e_{03}$&$-e_{02}$&$-e_{01}$&$e_{08}$&$-e_{07}$&$e_{06}$&$-e_{05}$\\
\hline
$e_{05}$&$e_{05}$&$-e_{06}$&$-e_{07}$&$-e_{08}$&$-e_{01}$&$e_{02}$&$e_{03}$&$e_{04}$\\
\hline
$e_{06}$&$e_{06}$&$e_{05}$&$-e_{08}$&$e_{07}$&$-e_{02}$&$-e_{01}$&$-e_{04}$&$e_{03}$\\
\hline
$e_{07}$&$e_{07}$&$e_{08}$&$e_{05}$&$-e_{06}$&$-e_{03}$&$e_{04}$&$-e_{01}$&$-e_{02}$\\
\hline
$e_{08}$&$e_{08}$&$-e_{07}$&$e_{06}$&$e_{05}$&$-e_{04}$&$-e_{03}$&$e_{02}$&$-e_{01}$\\
\hline
\end{tabular}\,\,,
$$
$$
\begin{tabular}{|r||r|r|r|r|r|r|r|r|}\hline
$\cdot$ &$e_{11}$&$e_{12}$&$e_{13}$&$e_{14}$&$e_{15}$&$e_{16}$&$e_{17}$&$e_{18}$\\
\hline\hline
$e_{01}$&$e_{11}$&$e_{12}$&$e_{13}$&$e_{14}$&$e_{15}$&$e_{16}$&$e_{17}$&$e_{18}$\\
\hline
$e_{02}$&$e_{12}$&$-e_{11}$&$e_{14}$&$-e_{13}$&$e_{16}$&$-e_{15}$&$-e_{18}$&$e_{17}$\\
\hline
$e_{03}$&$e_{13}$&$-e_{14}$&$-e_{11}$&$e_{12}$&$e_{17}$&$e_{18}$&$-e_{15}$&$-e_{16}$\\
\hline
$e_{04}$&$e_{14}$&$e_{13}$&$-e_{12}$&$-e_{11}$&$e_{18}$&$-e_{17}$&$e_{16}$&$-e_{15}$\\
\hline
$e_{05}$&$e_{15}$&$-e_{16}$&$-e_{17}$&$-e_{18}$&$-e_{11}$&$e_{12}$&$e_{13}$&$e_{14}$\\
\hline
$e_{06}$&$e_{16}$&$e_{15}$&$-e_{18}$&$e_{17}$&$-e_{12}$&$-e_{11}$&$-e_{14}$&$e_{13}$\\
\hline
$e_{07}$&$e_{17}$&$e_{18}$&$e_{15}$&$-e_{16}$&$-e_{13}$&$e_{14}$&$-e_{11}$&$-e_{12}$\\
\hline
$e_{08}$&$e_{18}$&$-e_{17}$&$e_{16}$&$e_{15}$&$-e_{14}$&$-e_{13}$&$e_{12}$&$-e_{11}$\\
\hline
\end{tabular}\,\,,
$$
$$
\begin{tabular}{|r||r|r|r|r|r|r|r|r|}\hline
$\cdot$ &$e_{01}$&$e_{02}$&$e_{03}$&$e_{04}$&$e_{05}$&$e_{06}$&$e_{07}$&$e_{08}$\\
\hline\hline
$e_{11}$&$e_{11}$&$e_{12}$&$e_{13}$&$e_{14}$&$\lambda e_{15}$&$\lambda e_{16}$&$\lambda e_{17}$&$\lambda e_{18}$\\
\hline
$e_{12}$&$e_{12}$&$-e_{11}$&$e_{14}$&$-e_{13}$&$\lambda e_{16}$&$-\lambda e_{15}$&$-\lambda e_{18}$&$\lambda e_{17}$\\
\hline
$e_{13}$&$e_{13}$&$-e_{14}$&$-e_{11}$&$e_{12}$&$\lambda e_{17}$&$\lambda e_{18}$&$-\lambda e_{15}$&$-\lambda e_{16}$\\
\hline
$e_{14}$&$e_{14}$&$e_{13}$&$-e_{12}$&$-e_{11}$&$\lambda e_{18}$&$-\lambda e_{17}$&$\lambda e_{16}$&$-\lambda e_{15}$\\
\hline
$e_{15}$&$e_{15}$&$-e_{16}$&$-e_{17}$&$-e_{18}$&$-\lambda e_{11}$&$\lambda e_{12}$&$\lambda e_{13}$&$\lambda e_{14}$\\
\hline
$e_{16}$&$e_{16}$&$e_{15}$&$-e_{18}$&$e_{17}$&$-\lambda e_{12}$&$-\lambda e_{11}$&$-\lambda e_{14}$&$\lambda e_{13}$\\
\hline
$e_{17}$&$e_{17}$&$e_{18}$&$e_{15}$&$-e_{16}$&$-\lambda e_{13}$&$\lambda e_{14}$&$-\lambda e_{11}$&$-\lambda e_{12}$\\
\hline
$e_{18}$&$e_{18}$&$-e_{17}$&$e_{16}$&$e_{15}$&$-\lambda e_{14}$&$-\lambda e_{13}$&$\lambda e_{12}$&$-\lambda e_{11}$\\
\hline
\end{tabular}
$$
and
$$ e_{1s}e_{1t}:=0 \quad\mbox{for $1\le s, \, t\le 8$}.$$

\hfill\raisebox{1mm}{\framebox[2mm]{}}

\medskip
$\mathcal{O}^{(2\lambda)}$ is both a non-associative composition $\mathcal{Z}_2$-algebra and an alternative algebra. The following is my conjecture about finite dimensional composition $\mathcal{Z}_2$-algebras.

\medskip
{\bf Conjecture} A finite dimensional composition $\mathcal{Z}_2$-algebra is isometrically isomorphic to one of the following ten $\mathcal{Z}_2$-algebras:
$$\mbox{
$\mathcal{R}$, \,\, $\mathcal{C}$, \,\, $\mathcal{H}$, \,\, $\mathcal{R}^{(2)}$,
\,\, $\mathcal{C}^{(2)}$, \,\, $\mathcal{C}^{(2)}$,  \,\, $\mathcal{H}^{(2)}$, 
\,\, $\mathcal{H}^{(-2)}$, \,\, $\mathcal{O}^{(2)}$ \, and \, $\mathcal{O}^{(-2)}$.}$$

\bigskip
\section{The Extended Lie Correspondence}

In this section we assume that all algebras are finite dimensional algebras over the field of real numbers.

\medskip
Let $G$ be a linear $\xi$-group in a Banach $\mathcal{Z}_2$-graded algebra $A$. The {\bf tangent space} $T_1(G)$ to $G$ at the identity $1$ of $A$ is defined by
$$
T_1(G):=\{\, a'(0) \, |\, \mbox{$a(t): \mathcal{R}\to G$ is a differential curve with $a(0)=1$ }\,\}.
$$ 
By Proposition 3.1 in \cite{hl1}, $T_1(G)$ is a Hu-Liu Leibniz algebra over the field of real numbers, which is also called the {\bf Hu-Liu Leibniz algebra of $G$}. 

\medskip
The next proposition gives an extension of the Campbell-Baker-Hausdorff series.

\begin{proposition}\label{pr3.1} Let $x$, $y$, $u$, $w$ be elements of a finite dimensional Banach $\mathcal{Z}_2$-algebra $A$. If $x$, $y$, $u$, $w$ are sufficiently close to $0$, then the equation 
$$
(exp \,u)_0(exp \,x)(exp \,u)_0^{-1}\cdot (exp \,w)_0(exp \,y)(exp \,w)_0^{-1}=exp \,z
$$
has a unique solution $z=C(x,\, y,\, u,\, w)$ in $A$ such that $z=C(x,\, y,\, 0,\, 0)$ is the Campbell-Baker-Hausdorff series.
\end{proposition}

\hfill\raisebox{1mm}{\framebox[2mm]{}}

\medskip
The first a few terms of the series $z=C(x,\, y,\, u,\, w)$ are as follows
\begin{eqnarray*}
&&C(x,\, y,\, u,\, w)=x+y+\frac12 [x, y]-\langle x, u\rangle-\langle y, w\rangle
 +\frac{1}{12}[x, [x, y]]+\\
&& +\frac{1}{12}[y, [y, x]]+\frac12 \langle\langle x, u\rangle , u\rangle
+\frac12 \langle\langle y, w\rangle , w\rangle
-\frac12 [x, \langle y, w\rangle ]-\frac12 [\langle x, u\rangle , y]+\\
&& +\frac{1}{12}[x, [x, y]]+\frac{1}{12}[y, [y, x]]
-\frac16 \langle \langle\langle x, u\rangle , u\rangle, u\rangle
-\frac16 \langle \langle\langle y, w\rangle , w\rangle, w\rangle +\\
&& +\frac14 [x,  \langle\langle y, w\rangle , w\rangle ]
+\frac12 [ \langle x, u\rangle , \langle y, w\rangle ]
+\frac14 [ \langle\langle x, u\rangle , u\rangle , y]
-\frac{1}{12}[x, [x, \langle y, w\rangle ] ]+\\
&& -\frac{1}{12}[ \langle x, u\rangle , [x, y] ]-\frac{1}{12}[y, [y, \langle x, u\rangle ] ]
-\frac{1}{12}[ \langle y, w\rangle , [y, x] ]+\\
&&+\frac{1}{24}[y, [x, [ y, x ] ] ]+\cdots 
\end{eqnarray*}

\medskip
We now state the extended Lie correspondence.

\begin{proposition}\label{pr3.2} Let $A$ be a finite dimensional Banach $\mathcal{Z}_2$-algebra, and let
$$
\mathcal{G}:=\{\, G\, |\, \mbox{$G$ is a connected linear $\xi$-group in $A$} \,\}
$$
and
$$
\mathcal{L}:=\{\, L\, |\, \mbox{$L$ is a linear Hu-Liu Leibniz algebra in $A$} \,\}.
$$
If $f: \mathcal{G}\to \mathcal{L}$ is the map defined by 
$$
f(G):=\mbox{the tangent space $T_1(G)$ to $G$ at the identity $1$ of $A$ for $G\in\mathcal{G}$},
$$
then $f$ is bijective, and the inverse map $f^{-1}: \mathcal{L}\to \mathcal{G}$ is given by
$$
f^{-1}(L):=\mbox{the linear $\xi$-group generated by $exp\, L$ for $L\in\mathcal{L}$}.
$$
\end{proposition}
\hfill\raisebox{1mm}{\framebox[2mm]{}}

\medskip
The extended Lie correspondence `$\xi$-groups $\leftrightarrow$ Hu-Liu Leibniz algebras' works not only on the level of sets, but also on the level of maps. In fact, both Theorem 2.21 and Theorem 3.7 in \cite{H} have extensions in the context of $\xi$-groups and Hu-Liu Leibniz algebras.

\end{document}